# Harmonic Center of a Closed Convex Polytope: Definition, Calculation and Some Properties

V. S. Patwardhan[1]

30 March 2020

**Abstract**

A closed convex polytope in n dimensions is defined by the usual linear inequality constraints specified as Ax ≤ b, where A is a (m x n) matrix, with m > n. If L is a straight line drawn in any direction from any feasible point P, then in general, it intersects every constraint at one point, giving m intersections. It is shown that there exists a unique feasible point Q somewhere along this line, such that $\sum_{i=1}^{i=m}(1/d_i) = 0$, where $d_i$ is the algebraic distance between Q and the intersection with the i[th] constraint, measured along the line. The point Q is defined as the harmonic point along the line L. The harmonic center of the polytope is defined as that point which is the harmonic point for all n lines drawn through it, each parallel to one of the coordinate axes. The existence and uniqueness of such a center is shown. The harmonic center can be calculated using the coordinate search algorithm (CS), as illustrated with some simple examples. The harmonic center defined here is a generalisation of the BI center defined earlier and is better in several respects. It is shown that the harmonic center of the polytope is also the harmonic point for any line drawn through it in any direction. It is also shown that for any strictly feasible point P, there exists a unique harmonic hyperplane passing through it, such that P is the harmonic point for any line which lies in the harmonic hyperplane and passes through P.

**1. Introduction**

The calculation of the center of a closed convex polytope has been of interest for a long time. One of the motivations for this interest has been the fact that interior point algorithms for convex programming, including linear programming, involve a centering step at periodic intervals. Therefore, any method which leads to efficient calculation of a uniquely defined center may become useful in interior point algorithms. Many different definitions have been used for the center, which include the center of mass of the entire polytope, the centroid of all vertices, the center of the largest inscribed sphere (or an ellipsoid), the center of the smallest sphere (or ellipsoid) which includes the polytope, the analytical center, the weighted projection center, orthogonal projections on to polytope faces and the BI center [Moretti, 2003; Patwardhan, 2019]. All these definitions lead to different points as centers, whose calculation involves different degrees of computational effort. Some of these, and a few others such as the Fermat-Torricelli point, orthocentre and the Monge point, have been investigated in detail for simplexes [Edmondsa et.al., 2004].

In this paper, we consider a closed convex polytope defined by the set of constraints Ax ≤ b, where A is a (m x n) matrix, with m > n. It is assumed here that all the m constraints are inequality constraints, i.e. there are no equality constraints. (Any equality constraint can be removed by eliminating any one variable using the equality constraint). We define a new center, i.e. the harmonic center, for a closed convex polytope. First, we define the harmonic point for an arbitrary line drawn through any feasible point within the polytope. The harmonic center of a polytope is then defined as that point which is the harmonic point for all n lines drawn through it, parallel to the coordinate axes. We illustrate a coordinate search algorithm for calculating the harmonic center, starting from any interior point. We show that the harmonic center is the harmonic point for any line drawn through it. The harmonic center defined here is a generalisation of the BI center defined earlier and is better in several respects.

---

[1] Independent researcher. Formerly, Scientist G, National Chemical Laboratory, Pune 411008, India.
   Email : vspatw@gmail.com , URL : https://www.vspatwardhan.com

## 2. Harmonic point along any line drawn through a feasible point

Let us illustrate the idea of a harmonic point along a line by considering a simple example shown in Figure 1 which shows a triangle (which is a simplex in two dimensions). Let us consider a strictly feasible point P which satisfies AP < b. Any line drawn through P, intersects the three constraints at points 1, 2 and 3. Q is a feasible point somewhere along this line. The distances between Q and points 1, 2 and 3 are shown as $d_1$, $d_2$ and $d_3$. It may be noted that $d_2$ is negative. The sum $\sum_{i=1}^{i=3}(1/d_i)$ goes from negative infinity to positive infinity as Q moves from 2 to 1 and is easily seen to be monotonic. Thus, there is a single feasible point along the line where this sum is equal to zero. This point (shown as Q in Figure 1) is defined as the harmonic point along the line.

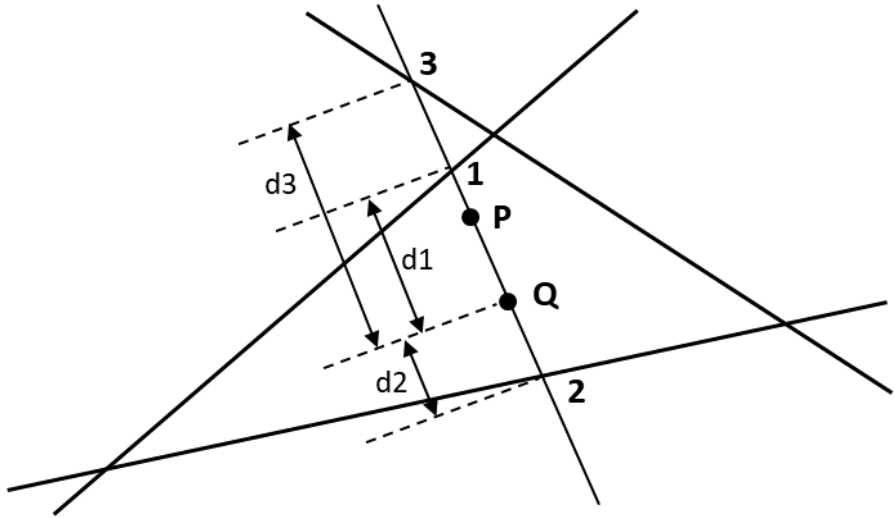

Figure 1. The harmonic point Q along a line through point P

## 3. Harmonic center of a closed convex polytope

The harmonic center of a polytope is defined as that feasible point which is the harmonic point for all n lines drawn through it, parallel to the coordinate axes. Let us illustrate the idea with an example in two dimensions. Figure 2 shows a triangle defined by three constraints. L is a line parallel to the $X_1$ axis and intersects the constraints at three points. The point Q is the harmonic point for this line. If the line L is moved vertically, Q shifts as well. The curve $C_1$ is the locus of Q. Any point on $C_1$ is the harmonic point for a line drawn through that point, parallel to the $X_1$ axis. Figure 3 shows two lines. The line $C_2$ is the locus of harmonic points for lines drawn parallel to the $X_2$ axis. Point H is the intersection of the two lines $C_1$ and $C_2$ and is the harmonic point for both lines drawn through it, parallel to the two axes. Thus, H is the harmonic center of the triangle.

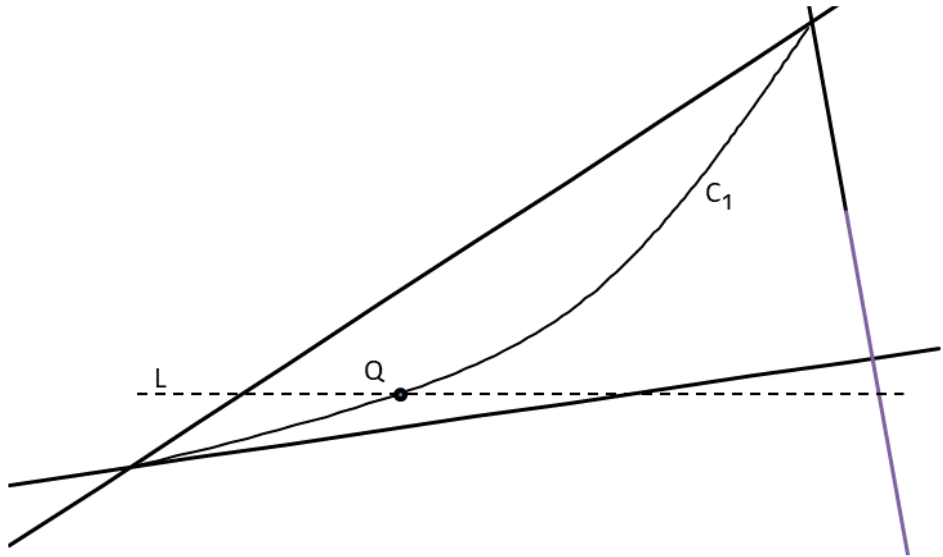

Figure 2. Locus $C_1$ of harmonic points along horizontal lines

The existence and uniqueness of the harmonic center for a closed convex polytope can be shown using arguments like those used earlier by Patwardhan [2] in relation to the bisection center.

**4. Calculation of the harmonic point along a line parallel to a coordinate axis**

For a closed convex polytope in n dimensions defined by m inequality constraints, the i$^{th}$ constraint can be written as

$$\sum_{j=1}^{n} A_{ij} x_j \leq b_i \qquad for\ i = 1\ to\ m \tag{1}$$

It is assumed here that the rows of A are normalized. In other words, $\sum_{j=1}^{j=n} A_{ij}^2 = 1\ \ for\ all\ i$. Let P be a strictly feasible interior point so that we can write

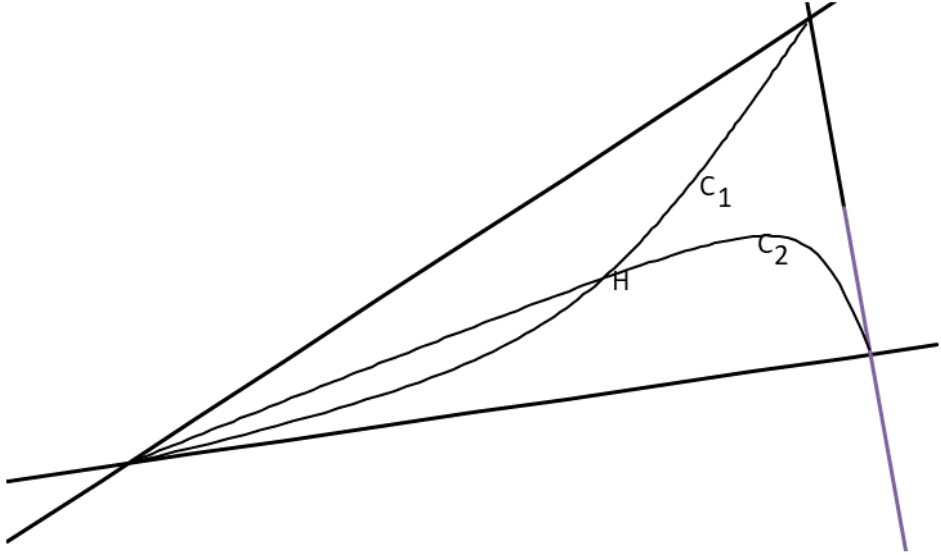

Figure 3. Harmonic center as intersection of lines $C_1$ and $C_2$

$$\sum_{j=1}^{n} A_{ij}P_j < b_i \quad for\ i = 1\ to\ m \tag{2}$$

The residuals at point P can be denoted as $S_i$, which are given by

$$S_i = b_i - \sum_{j=1}^{n} A_{ij}P_j \tag{3}$$

Since P is an interior point, $S_i > 0$ for all i. A line drawn through P, parallel to the $k^{th}$ axis, intersects the $i^{th}$ constraint after travelling a distance $d_{ik}$ given by

$$d_{ik} = S_i/A_{ik} \quad for\ k = 1\ to\ n \tag{4}$$

Since the rows of A are assumed to be normalized, $d_{ik}$ is the actual distance. If $A_{ik} = 0$, then $d_{ik}$ becomes infinite, i.e. there is no intersection with constraint i. If $A_{ik} > 0$, then $d_{ik} > 0$, i.e. the line drawn in the $+x_k$ direction intersects the $i^{th}$ constraint. On the other hand, if $A_{ik} < 0$, then $d_{ik} < 0$, i.e. the line drawn in the $-x_k$ direction intersects the $i^{th}$ constraint. Let $d_{+k}$ and $d_{-k}$ be the distances to the first intersection in the $+x_k$ and $-x_k$ directions respectively, and let $Q_{+k}$ and $Q_{-k}$ be the corresponding points of intersection. These are given by

$$d_{+k} = \min(d_{ik}) \quad over\ all\ i, for\ which\ A_{ik} > 0 \tag{5}$$

and

$$d_{-k} = \max(d_{ik}) \quad \text{over all } i, \text{for which } A_{ik} < 0 \tag{6}$$

It is obvious that $d_{+k} > 0$, and $d_{-k} < 0$.

Let Q be a point on this line which is the harmonic point for this line. Let $Q_k = P_k + h_k$. Then, using the definition of a harmonic point, we can write

$$F_k = \sum_{i=1}^{m} \frac{1}{(d_{ik} - h_k)} = 0 \tag{7}$$

$F_k$ becomes minus infinity at $h_k = d_{-k}$, and plus infinity at $h_k = d_{+k}$. $F_k$ is zero only for a single value of $h_k$. This equation can be easily solved for $h_k$ using the Newton Raphson method, and converges within a few iterations, giving the point Q. The procedure can be repeated for lines parallel to other coordinate axes.

## 5. Harmonic center of a polytope

It has been mentioned before that the harmonic center of a closed convex polytope is defined as that feasible point which is the harmonic point for all n lines drawn through it, parallel to the coordinate axes. Thus, the harmonic center H is characterised by

$$F_k = \sum_{i=1}^{m} \frac{1}{d_{ik}} = 0 \qquad \text{for } k = 1 \text{ to } n \tag{8}$$

where $d_{ik}$ is given by equation (4).

The harmonic center has an interesting property. Let us consider any arbitrary line with the unit vector $U = [u_1, u_2, \ldots u_n]^T$. The distance from H to the constraint i along this line is given by

$$d_i = \frac{S_i}{\sum_{j=1}^{n} A_{ij} u_j} \tag{9}$$

Let us consider the sum

$$F_u = \sum_{i=1}^{m} \frac{1}{d_i} = \sum_{i=1}^{m} \left[ \frac{\sum_{j=1}^{n} A_{ij} u_j}{S_i} \right] \tag{10}$$

After some algebraic rearrangement, this can be written as

$$F_u = \sum_{i=1}^{m} \left[ \frac{\sum_{j=1}^{n} A_{ij} u_j}{S_i} \right] = \sum_{j=1}^{n} \left[ u_j \sum_{i=1}^{m} \frac{A_{ij}}{S_i} \right] = \sum_{j=1}^{n} u_j F_j \tag{11}$$

Since H is the harmonic center, $F_j = 0$ for all j, as shown earlier. Therefore, $F_u = 0$. Thus, the harmonic center H of the polytope is also the harmonic point for any arbitrary line drawn through H.

The harmonic center of a closed convex polytope has been defined earlier as that feasible point which is the harmonic point for all n lines drawn through it, parallel to the coordinate axes. In view of this result, it can be defined more simply as that feasible point which is the harmonic point for any arbitrary line drawn through it.

## 6. Harmonic hyperplane for a given feasible point

Closed convex polytopes have another interesting property. It is shown below that there exists a unique hyperplane corresponding to any given feasible point P, such that P is the harmonic point for any arbitrary line which passes through P and lies in the harmonic hyperplane.

Let P be any strictly feasible point. Let us consider any arbitrary line with the unit vector U = [$u_1$, $u_2$, …. $u_n$]$^T$. The distance from P to the constraint i along this line is given by

$$d_i = \frac{S_i}{\sum_{j=1}^{n} A_{ij} u_j} \quad (12)$$

where $S_i$ is given by Equation (3). Let us consider the sum

$$F_u = \sum_{i=1}^{m} \frac{1}{d_i} = \sum_{i=1}^{m} \left[\frac{\sum_{j=1}^{n} A_{ij} u_j}{S_i}\right] \quad (13)$$

After some algebraic rearrangement, this can be written as

$$F_u = \sum_{i=1}^{m} \left[\frac{\sum_{j=1}^{n} A_{ij} u_j}{S_i}\right] = \sum_{j=1}^{n} \left[u_j \sum_{i=1}^{m} \frac{A_{ij}}{S_i}\right] = \sum_{j=1}^{n} u_j v_j \quad (14)$$

where

$$v_j = \sum_{i=1}^{m} \frac{A_{ij}}{S_i} \quad for\ j = 1,2,…n \quad (15)$$

It may be noted that vector V, with components given by Equation (15), is uniquely defined by the point P. We want P to be the harmonic point for the line U passing through it. Therefore, $F_u$ = 0, and vector U must be orthogonal to vector V, as implied by Equation (14). Thus, the vector U must lie in the hyperplane defined by

$$\sum_{j=1}^{n} v_j x_j = \sum_{j=1}^{n} v_j P_j \quad (16)$$

which can be termed as the harmonic hyperplane for P. This proves that there exists a unique hyperplane corresponding to any given feasible point P (as given by Equation (16)), such that P is the harmonic point for any arbitrary line U which passes through P and lies in the harmonic hyperplane.

## 7. Characterisation of the harmonic center

For a strictly feasible point P, we define a vector F whose components are given by

$$F_j = \sum_{i=1}^{m} \frac{1}{d_{ij}} = 0 \qquad for\ j = 1\ to\ n \qquad (17)$$

The magnitude of F, i.e. $\|F\|$, can be used as a good indicator of the closeness of an interior point P to the harmonic center H. It may be noted that $\|F\|$ = 0 at the harmonic center, ∞ at the polytope boundary, and ≥ 0 everywhere else in the feasible region.

The problem of finding the harmonic center is essentially one of minimising $\|F\|$ (to a value of 0) using unconstrained optimisation. It may be noted that $\|F\|$ is continuous and has continuous derivatives over the entire feasible region. Therefore, optimisation techniques based on using the first derivative can, in principle, be used. Unfortunately, analytical derivatives are algebraically too complex, and even numerical derivatives involve substantial computational effort. We take a simplified approach for determining the harmonic center, which uses a coordinate search algorithm. A good recent introduction to coordinate search is available [Shi et al., 2017].

## 8. Coordinate search algorithm (CS)

The idea here is to go from a given interior point P to a new point Q in n stages, changing one coordinate at a time to go to the harmonic point along that coordinate axis, and repeat this step iteratively till we approach the harmonic center sufficiently closely, as indicated by $\|F\|$. Several variations of the CS algorithm are described by Astolfi [2006].

One step of this algorithm consists of the following: Given an interior point P, we get the point $Q_1$, which is the harmonic point of the line through P in the $x_1$ direction, which reduces $\|F\|$. Then we get the point $Q_2$ which is the harmonic point of the line through $Q_1$ in the $x_2$ direction, which reduces $\|F\|$ further. We continue this till we get point $Q_n$, which is the harmonic point of the line through $Q_{n-1}$ in the $x_n$ direction. This completes the step. It may be noted that at each stage we go to a harmonic point of some line, thereby ensuring that we do not approach the polytope boundary too closely. This process is continued till we approach the harmonic center closely. The details of the algorithm are given below.

___________________________________________________________

CS Algorithm:
    *Given : iter = 0, and $P_{iter}$ = $P_0$ (an interior point)*
     *Do while $\|F\|$ at $P_{iter}$ > a critical value (such as 0.01)*
      *Q(j) = $P_{iter}$(j) for j = 1 to n*
       *For k = 1 to n*
        *Q(k) = Q(k) + $h_k$ with $h_k$ from eq. (7)*
        *Update $S_i$ at Q(k) for all i*
       *Next k*
      *iter = iter + 1*
      *$P_{iter}$(j) = Q(j) for j = 1 to n*
     *End do*

*P<sub>iter</sub> is the harmonic center of the polytope*

___

Let us illustrate the progress of this algorithm with following two examples.

___

Example 1:   A polytope with 5 constraints and 2 nonnegative variables

$$\begin{aligned} 1.5\,x - y &\leq 8 \\ 0.2\,x + y &\leq 8.4 \\ -5\,x - y &\leq -10 \\ -4\,x + y &\leq 1 \\ 0.5\,x - y &\leq 2 \\ x,\ y &\geq 0 \end{aligned} \qquad (18)$$

___

Let us use six starting points, i.e. (9, 6), (3, 0.25), (2, 7.5), (5, 1), (7, 3) and (5, 7). Figure 4 shows the progress of this algorithm for the six starting points. All the lines converge to the harmonic center H, i.e. (6.02, 5.55). Table 1 shows the results of iterative calculations for all these cases. It is seen that the calculations converge to the harmonic center in just two iterations.

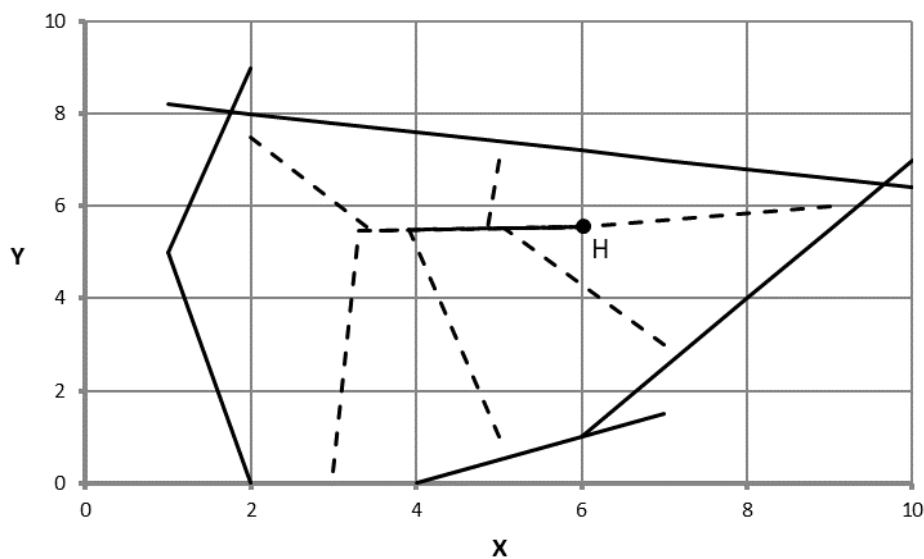

Figure 4. Approach to harmonic center from several starting points

| No. | Starting point | | Iteration 1 | | Iteration 2 | |
|---|---|---|---|---|---|---|
| | x | y | x | y | x | y |
| 1 | 9.00 | 6.00 | 6.01 | 5.55 | | |
| 2 | 3.00 | 0.25 | 3.31 | 5.47 | 6.02 | 5.55 |
| 3 | 2.00 | 7.50 | 3.46 | 5.48 | 6.02 | 5.55 |
| 4 | 5.00 | 1.00 | 3.91 | 5.49 | 6.02 | 5.55 |
| 5 | 7.00 | 3.00 | 5.07 | 5.52 | 6.03 | 5.55 |
| 6 | 5.00 | 7.00 | 4.86 | 5.51 | 6.03 | 5.55 |

Table 1: Approach to the harmonic center for Example 1

Let us consider a larger problem in four dimensions as given below.

______________________________________________________

Example 2:    A polytope with 5 constraints and 4 nonnegative variables

$$\begin{aligned} x_1 + x_2 - x_3 + x_4 &\leq 8 \\ x_1 + 0.5\, x_2 - x_3 &\leq 3 \\ 0.5\, x_1 - 2\, x_2 + x_3 &\leq 2 \\ - x_1 + 0.5\, x_2 - 0.5\, x_4 &\leq 3 \\ x_1 + 3\, x_2 + 1.5\, x_3 + 2\, x_4 &\leq 25 \\ x_1,\ x_2,\ x_3,\ x_4 &\geq 0 \end{aligned}$$

(19)

______________________________________________________

This example has five constraints defining the polytope, and in addition, all four variables are nonnegative, which is equivalent to four more constraints.

If we apply the CS algorithm to this polytope, starting from the arbitrary point (1, 2, 2.5, 1.3) we approach the harmonic center iteratively. Figure 5 shows the variation of the four coordinates with iteration number. Table 2 shows the approach to the harmonic center with iteration number. It also shows $\|F\|$, which indicates the closeness of the point to the harmonic center. It is seen that the approach to the harmonic center is fast, and it takes only three iterations.

## 9. Comparison with the BI center

The BI center of a closed convex polytope [Patwardhan, 2019] was defined as the midpoint of $Q_{+k}$ and $Q_{-k}$, calculated using Equations (5) and (6), for all directions parallel to the coordinate axes. These two points represent the nearest intersections in the positive and negative directions. The harmonic center defined here is based on all the intersections. It was shown earlier that the function which gets maximised at the BI center, is continuous but only piecewise differentiable. This can lead to difficulties in convergence. The harmonic center defined here minimises $\|F\|$ to zero, which is continuous as well as differentiable, and gives a smoother approach to the minimum. The BI center does not bisect all lines passing through it. It only bisects lines drawn through it in directions parallel to the coordinate axes. The harmonic center, on the other hand, is a harmonic point for any line passing through it, and not just for lines parallel to the coordinate axes. The harmonic center is thus a generalisation of the BI center proposed earlier and is expected to give computational advantages for polytopes in higher dimensions.

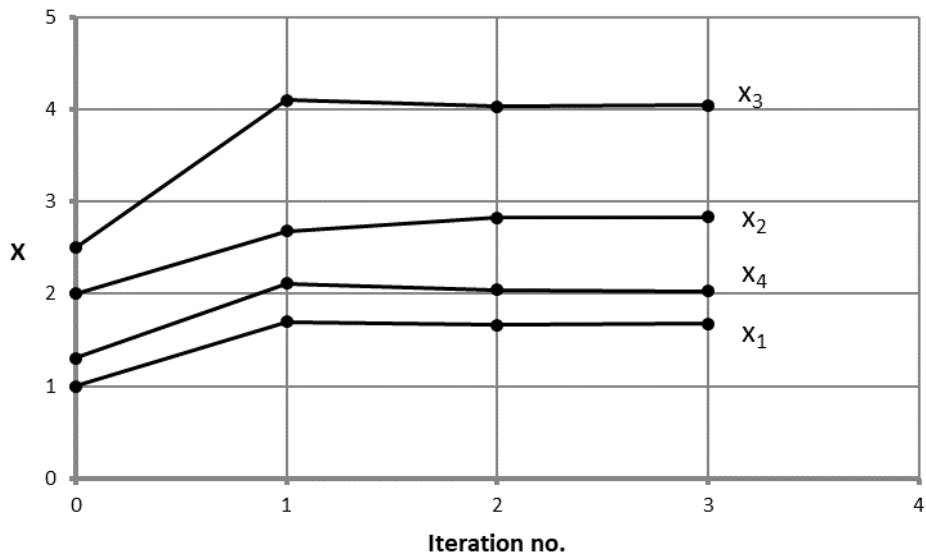

Figure 5. Variation of x values with iteration number

| Iter. No. | $x_1$ | $x_2$ | $x_3$ | $x_4$ | $\|F\|$ |
|---|---|---|---|---|---|
| 0 | 1.00 | 2.00 | 2.50 | 1.30 | 1.013 |
| 1 | 1.70 | 2.68 | 4.10 | 2.11 | 0.185 |
| 2 | 1.66 | 2.82 | 4.03 | 2.05 | 0.016 |
| 3 | 1.68 | 2.83 | 4.05 | 2.03 | 0.007 |

Table 2: variation of X values and $\|F\|$ with iteration number

## 10. Conclusions

A new center, termed here as the harmonic center, has been defined for a closed convex polytope. The harmonic center is defined as that point which is the harmonic point for all n lines drawn through it, each parallel to one of the coordinate axes. The calculation of the harmonic center has been illustrated using the coordinate search algorithm, for some simple examples. The harmonic center defined here is a generalisation of the BI center defined earlier and is better in several respects. The harmonic center is also the harmonic point for any line drawn through it in any direction. For any interior point P, there exists a unique harmonic hyperplane passing through it, such that P is the harmonic point for any line which lies in the harmonic hyperplane and passes through P.

In this paper, calculation of the harmonic center has been illustrated using the coordinate search algorithm for two small problems. These ideas need to be applied to much larger problems for developing numerically efficient procedures. It is likely that the properties of a harmonic center described above can be utilised for this purpose. This work is in progress at present.

## 11. References

[1] A. C. Moretti,
*A weighted projection centering method*,
Computational and Applied Mathematics, **22** (2003), pp 19–36

[2] Patwardhan V. S.,
*Center Of A Closed Convex Polytope: A New Definition*,
DOI: 10.13140/RG.2.2.35394.53446 (2019)
(Can also be downloaded from www.vspatwardhan.com )

[3] Patwardhan V. S.,
*Calculation of the BI Center of a Closed Convex Polytope*,
DOI: 10.13140/RG.2.2.32331.95523 (2019)
(Can also be downloaded from www.vspatwardhan.com )

[4] A. L. Edmondsa, M. Hajjab and H. Martin,
*Coincidences of simplex centers and related facial structures*,
http://arxiv.org/abs/math/0411093v1 (2004)

[5] Shi H. M., Tu S., Xu Y. and Yin W.,
*A Primer on Coordinate Descent Methods*,
arXiv:1610.00040 [math.OC] , (2017)

[6] Astolfi A.,
*Optimisation: An Introduction,*
http://www3.imperial.ac.uk/pls/portallive/docs/1/7288263.PDF , (2006)